\documentclass[reqno]{amsart}
\usepackage{amssymb}
\usepackage{hyperref}
\newtheorem{definition}{Definition}[section]
\newtheorem{thm}{Theorem}[section]
\newtheorem{proposition}{Proposition}[section]
\newtheorem{corollary}{Corollary}[section]
\newtheorem{lemma}{Lemma}[section]

\newtheorem{remark}{Remark}[section]

\begin{document}

\title[Some generalized inequalities]{Some generalized inequalities involving extended beta and gamma functions for several variables}

\author[S. Mubeen, I. Aslam, Ghazi S. Khammash, Saralees Nadarajah, Ayman Shehata]{S. Mubeen, I. Aslam, Ghazi S. Khammash, Saralees Nadarajah, Ayman Shehata}

\address{S. Mubeen, Department of Mathematics, University of Sargodha, Sargodha, Pakistan}
\email{smjhanda@gmail.com}
\address{I. Aslam, Department of Mathematics, University of Sargodha, Sargodha, Pakistan}
\email{iaslam@gmail.com}
\address{Ghazi S. Khammash,  Department of Mathematics, Al-Aqsa University, Gaza Strip, Palestine}
\email{ghazikhamash@yahoo.com}
\address{Saralees Nadarajah (corresponding author), Department of Mathematics, University of Manchester, Manchester M13 9PL, UK}
\email{mbbsssn2@manchester.ac.uk}
\address{Ayman Shehata,   Department of Mathematics, Faculty of Science, Assiut University, Assiut 71516, Egypt}
\email{aymanshehata@science.aun.edu.eg, drshehata2006@yahoo.com}

\date{}
\maketitle{}

\begin{abstract}
Recently,   extensions of gamma and beta functions have been studied by
many  researchers due to their nice properties and variety of applications in different fields of science.
The aim of this note is to investigate  generalized inequalities associated with extended beta and gamma functions.
\end{abstract}



\section{Introduction}

In many branches of both pure and practical mathematics,
different kinds of special functions have been   developed as crucial tools for scientists and engineers.
Gamma and beta functions are the most notable  special functions,
Swiss mathematician Leonhard Euler initially introduced the gamma function.
It also arises in numerous other contexts like Riemann's zeta function,
asymptotic series, definite integration, hypergeometric series, and number theory.
Due to its significance, several renowned mathematicians like Adrien-Marie Legendre (1752-1733),
Carl Friedrich Gauss (1777-1755), Christoph Gudermann (1798-1752),
and  Joseph Liouville (1809-1782) have investigated the gamma function.
The special transcendental functions category include the gamma function.

The beta function has applications in many branches of science
and mathematics and was first researched by Euler and Legendre.
Jacques Binet gave it its name.
The function computes, for instance, the scattering amplitudes
of the Regge trajectories in string theory, a branch of complicated physics.

The fascinating and important area of mathematics known as special functions
has many uses in a variety of industries such as records, engineering, astronomy, etc.
Many projects on this topic have already been completed.
Numerous researchers are working on special functions and summation theorems.

Due to the variety of applications of beta and gamma functions
many researchers have derived their representations and properties.
Diaz \emph{et al.} \cite{m1}-\cite{m3} provided integral representations of beta and gamma $k$  functions
and derived their properties.
They also provided a  representation for the Pochhammer's $k$ symbol.
After the result,  many other researchers including
Kokologiannaki \cite{m4}-\cite{m5}, Krasiniqi \cite{m6}, Mansour \cite{m7}
and Mubeen \emph{et al.} \cite{m8} added their contributions, making these functions more interesting and useful.

Mubeen \emph{et al.} \cite{m9} discussed a representation of the
beta and gamma $k$  functions.
Golub \cite{m10}  contributed in this frame work.
Mubeen and Habibullah \cite{m11} provided  integral representations
of some generalized confluent hypergeometric $k$  function by using properties of Pochhammer's $k$  symbol,
beta function   and gamma $k$  function.
Mubeen \emph{et al.} \cite{m12} studied  other extensions of gamma and beta $k$  functions
involving a  confluent hypergeometric $k$ function.
Mubeen \cite{m13} introduced a  $k$  analogue of Kommer's formula and evaluated some useful result
by using hypergeometric $k$  functions.

Rehman \emph{et al.} \cite{m14}  introduced a  beta $k$ function for several variables.
They also extended beta $k$ function for $n$ variables.
Mubeen and Habibullah \cite{m15} defined the $k$ fractional integration and gave its application.
Rehman \emph{et al.} \cite{m16}-\cite{m17} derived some inequalities involving beta and gamma $k$  functions.
Raissouli and Soubhy \cite{m18} studied some inequalities involving two generalized beta functions in $n$ variables.

The beta and gamma $k$ functions have the following standard representations  \cite{m9}
\begin{equation*}
\displaystyle
\beta_{k}{(\phi,\psi)}=\frac {1}{k}\int\limits_{0}^{1}m^{\frac {\phi}{k}-1}{(1-m)}^{\frac {\psi}{k}-1}dm,
\qquad
where
\qquad
Re(\phi)>0, \qquad Re(\psi)>0
\end{equation*}
and
\begin{equation*}
\displaystyle
\Gamma_{k}{(\phi)}=\int\limits_{0}^{\infty}m^{{\phi}-1}e^{-\frac {m^{k}}{k}}dm,
\qquad
\textnormal {where}
\quad
Re(\phi)>0,
\qquad
k>0.
\end{equation*}
Some other extensions of extended beta and gamma $k$ functions  described in \cite{m4}-\cite{m8} are
\begin{equation*}
\displaystyle
\beta_{k}{(\phi,\psi;a)}=\frac {1}{k}\int\limits_{0}^{1}m^{\frac {\phi}{k}-1}{(1-m)}^{\frac {\psi}{k}-1}e^{-\frac {a^{k}}{k m(1-m)}}dm,
\qquad
where
\qquad
Re(\phi)>0,
\qquad
Re(\psi)>0,
\end{equation*}
\begin{equation}
\label{1.4}
\displaystyle
\beta_{k}{(\phi,\psi;a,b)}=\frac {1}{k}\int\limits_{0}^{1}m^{\frac {\phi}{k}-1}e^\frac {-a^{k}}{m}{(1-m)}^{\frac {\psi}{k}-1}e^{-\frac {(b)^{k}}{k(1-m)}}dm,
\qquad
where
\qquad
Re(\phi)>0,
\qquad
Re(\psi)>0
\end{equation}
and
\begin{equation}
\label{1.3}
\displaystyle
\Gamma_{a,k}{(\phi)}=\int\limits_{0}^{\infty}m^{{\phi}-1}e^{-\frac {m^{k}}{k}}e^{-\frac {a^{k}}{km^{k}}}dm,
\qquad
\textnormal {where}
\quad
Re(\phi)>0.
\end{equation}
The generalized beta and gamma functions  can be expressed as \cite{m12}
\begin{equation}
\label{1.5}
\displaystyle
\beta_{a,k}^{\left(a_{n},b_{n}\right)}{(\phi,\psi)}=\frac {1}{k}\int\limits_{0}^{1}m^{\frac {\phi}{k}-1}{(1-m)}^{\frac {\psi}{k}-1}
\ {}_{1}F_{1,k}\left(a_{n},b_{n};-\frac {a^{k}}{km(1-m)}\right) dm
\end{equation}
and
\begin{equation}
\label{1.6}
\displaystyle
\Gamma_{k}^{\left(a_{n},b_{n}\right)}{(\phi,a)}=\int\limits_{o}^{\infty}m^{\phi-1}
\ {}_{1}F_{1,k}\left(a_{n};b_{n};-\frac {m^{k}}{k}-\frac {a^{k}}{{km^{k}}}\right)dm,
\qquad
\textnormal
where
\quad a,b,>0.
\end{equation}
where ${}_{1}F_{1,k}$ is  the confluent hypergeometric function  \cite{m11} defined by
\begin{equation}
\label{1.7}
\displaystyle
{}_{1}F_{1,k} \left(a_{n},b_{n};l\right)=\sum_{m=0}^{\infty}\frac {\left(a_{n}\right)_{m,k}}{\left(b_{n}\right)_{m,k}}\cdot\frac {l^{m}}{m!}.
\end{equation}
If $a_{n} > 0$, $b_{n} - a_{n} >0$ and $k > 0$ then we have following integral representation
\begin{equation}
\label{1.8}
\displaystyle
{}_{1}F_{1,k} \left(a_{n},b_{n};l\right)=\frac {1}{k}
\frac {\Gamma\left(b_{n}\right)}{\Gamma \left(a_{n}\right) \Gamma\left(b_{n}-a_{n}\right)}
\int\limits_{0}^{1}u^{a_{n}-1}(1-u)^{b_{n}-a_{n}-1}e^{lu} du.
\end{equation}
By making the substitution, we obtain
\begin{equation}
\label{1.9}
\displaystyle
{}_{1}F_{1,k} \left(a_{n};b_{n};l\right)=e^{l} \ {}_{1}F_{1,k} \left(b_{n}-a_{n};b;-l\right).
\end{equation}

\begin{remark}
\label{11}
The real map is strictly increasing and strictly convex on $\Re$.
It follows that ${}_{1}F_{1,k} \left(a_{n};b_{n};l\right) \geq  {}_{1}F_{1,k} \left(a_{n};b_{n};0\right)$
for any $l\geq0$  and $0\leq  \  {}_{1}F_{1,k} \left(a_{n};b_{n};l\right)\leq 1$ for any $l\leq0$.
\end{remark}

\cite{m14} introduced the beta $k$ function for more then two variables
and provided some useful representations.
Let  $n > 3$ be an integer and let  $E_{(n-1)}$
be the $(n-1)$ simplex of $\Re^{n-1}$ described as
\begin{equation*}
\displaystyle
E_{(n-1)} = \left[ \left(m_{1}, \ldots,  m_{(n-1)} \right) \in
\Re^{(n-1)}:\sum_{i=1}^{n-1}m_{i}\leq 1;m_{i}\geq0;
\mbox{ for }
i=1, \ldots, n-1 \right].
\end{equation*}
The beta function involving $n$
variables $\phi_{1}, \ldots, \phi_{n}>0$ is
\begin{equation}
\label{1.10}
\displaystyle
\beta_{k} \left(\phi_{1}, \ldots, \phi_{n} \right) =
\frac {1}{k^{n-1}}\int_{\left(E_{n-1}\right)}
\prod_{i=1}^{n}m_{i}^{\frac {\phi_{i}}{k}}dm_{1} \cdots  dm_{n-1}.
\end{equation}
Let
\begin{equation*}
\displaystyle
m_{n}=1-\sum_{i=1}^{n-1}m_{i}
\end{equation*}
and
\begin{equation*}
\displaystyle
\sigma{(l)}=\sum_{i=1}^{n}{\left(\phi_{i}\right)},
\end{equation*}
then (\ref{1.10}) can be written as
\begin{equation*}
\displaystyle
\beta_{k} \left(\phi_{1}, \ldots, \phi_{n}\right)=\frac {\prod_{i=1}^{n}\Gamma_{k}\left(\phi_{i}\right)}{\Gamma_{k}\left(\sigma(\phi)\right)}.
\end{equation*}
Some other representations of the  extended beta $k$ function are
\begin{equation}
\label{1.14}
\displaystyle
\beta_{k} \left(\phi_{1}, \ldots,  \phi_{n};a\right) =
\frac {1}{k^{n-1}}\int\limits_{\left(E_{N-1}\right)}
\prod_{i=1}^{n}m_{i}^{\frac {\phi_{i}}{k}-1}e^{\frac {-a^{k}}{k\pi{(m)}}}dm_{1} \cdots  dm_{n-1}
\end{equation}
and
\begin{equation}
\label{1.15}
\displaystyle
\beta_{k} \left(\phi_{1}, \ldots,  \phi_{n};a_{1}, \ldots,  a_{n} \right) =
\frac {1}{k^{n-1}}\int\limits_{\left(E_{N-1}\right)}
\prod_{i=1}^{n}m_{i}^{\frac {\phi_{i}}{k}-1}e^{\frac {-a_{i}^{k}}
{k\pi(m)}}dm_{1} \cdots   dm_{n-1}
\end{equation}
for any $\phi_{1}, \ldots,  \phi_{n}>0$.

Now we defined the gamma $k$  function for the several variables.
Let $\phi =: \left(\phi_{1}, \ldots, \phi_{n}\right) >  0$,
$\alpha =: \left(\alpha_{1}, \ldots,  \alpha_{n} \right) > 0$,
$\beta=: \left(\beta_{1}, \ldots,  \beta_{n} \right)$
and $c = \left( c_{1}, \ldots,  c_{n} \right) \geq  0$.
The generalized gamma $k$  function is defined by
\begin{equation}
\label{2.1}
\displaystyle
\Gamma_{k,c}{(\phi)}=\int\limits_{(0,\infty)^{n}}
\prod_{i=1}^{n}m_{i}^{\phi_{i}-1}e^{-\frac {m_{i}^{k}}{k}}e^{-\frac {c_{i}^{k}}{km_{i}^{k}}}dm.
\end{equation}
If $ c=0$ then
\begin{equation*}
\displaystyle
\Gamma_{k, 0}{(\phi)}=\prod_{i=1}^{n}\Gamma_{k} \left(\phi_{i}\right).
\end{equation*}
Another representation of  (\ref{2.1}) is
\begin{equation}
\label{2.3}
\displaystyle
\Gamma_{k,c}^{\left(\alpha_{n},\beta_{n}\right)}{(\phi)} =
\int\limits_{(0,\infty)^{n}}\prod_{i=1}^{n}m_{i}^{{\phi_{i}}-1}
\ {}_{1}F_{1,k} \left( \left(\alpha_{n}\right)_{i}, \left(\beta_{n}\right)_{i}, -\frac {m_{i}^{k}}{k}-\frac {c_{i}^{k}}{km_{i}^{k}} \right) dm.
\end{equation}
If $c = 0$ then
\begin{equation}
\label{2.4}
\displaystyle
\Gamma_{k, 0}^{\left(\alpha_{n},\beta_{n}\right)}{(l)}=
\int\limits_{(0,\infty)^{n}}\prod_{i=1}^{n}m_{i}^{{\phi_{i}}-1}
\ {}_{1}F_{1}\left( \left(\alpha_{n}\right)_{i}; \left(\beta_{i}\right)_{n};-\frac {m_{i}^{k}}{k}\right) dm.
\end{equation}
If $n=1$ then (\ref{2.1}) and (\ref{2.3})  reduce to (\ref{1.3}) and (\ref{1.6}), respectively.
If $\alpha_{n} = \beta_{n}$ then (\ref{2.3}) becomes (\ref{2.1}).

\section{Main results}

In this section, we study the generalized beta $k$ function of the first kind.

\subsection{Generalized beta $k$  function of the first kind}

\begin{definition}
Let $\phi = \left(\phi_{1}, \ldots, \phi_{n} \right) > 0$, $a_{n} > 0$, $b_{n} > 0$, $\eta > 0 $ and $\zeta\geq0$.
The  generalized beta $k$ function of first kind is
\begin{equation}
\label{3.1}
\displaystyle
\beta_{\zeta,k}^{\left(a_{n},b_{n}\right)}{(\phi;\zeta)}=:
\frac {1}{k^{n-1}}\int\limits_{E_{n-1}}\prod_{i=1}^{n}t_{i}^{\frac {\phi_{i}}{k}-1}
\ {}_{1}F_{1,k}\left(a_{n},b_{n};-\frac {\eta^{k}}{k\pi(t)}-\zeta^{k}\frac {\pi(t)}{\eta^{k}}\right)dt,
\end{equation}
where  $dt =: dt_{1} \cdots  dt_(n-1)$ and $\pi(t)=:\prod_{i=1}^{n}t_{i}$ with $t_{n}=1-\sum_{i=1}^{n-1}t_{i}$.
If $\zeta = 0$, we obtain
\begin{equation}
\label{3.2}
\displaystyle
\beta_{0,  k}^{\left(a_{n},b_{n}\right)}{(\phi,\eta)} =
\frac {1}{k^{n-1}}\int\limits_{\left(E_n-1\right)}\prod_{i=1}^{n}t_{i}^{\frac {\phi_{i}}{k}-1}
\  {}_{1}F_{1,k} \left(a_{n},b_{n};-\frac {\eta^{k}}{k\pi(t)}\right)dt.
\end{equation}
If   $n=2$ and $\zeta=0$ then (\ref{3.1}) becomes (\ref{1.5}).
If $a_{n} = b_{n}$ then (\ref{3.1}) is (\ref{1.14}).
\end{definition}

\begin{proposition}
Let $\phi=\left(\phi_{1}, \ldots, \phi_{n}\right) > 0$, $a_{n} > 0$, $b_{n} > 0$,
$\eta > 0$  and $\zeta \geq 0$  with $ b_{n}-a_{n} > 0$.
Then, $ 0 \leq \beta_{\zeta,k}^{\left(a_{n},b_{n}\right)} \leq \beta_{k}{(\phi)}$ and so,
$\beta_{\eta,k}^{\left(a_{n},b_{n}\right)}{(\phi;\eta)}$ is well defined.
\end{proposition}

\begin{proof}
With the help of  Remark \ref{11}, we have
\begin{equation}
\label{3.3}
\displaystyle
0\leq\prod_{i=n}^{n}t_{i}^{\frac {\phi_{i}}{k}-1}
\ {}_{1}F_{1,k} \left(a_{n},b_{n};-\frac {\eta^{k}}{k\pi(t)}-\zeta^{k}\frac {k\pi(t)}{\eta^{k}}\right) \leq \
\prod_{i=1}^{n}t_{i}^{\frac {\phi_{i}}{k}-1}.
\end{equation}
Integrating (\ref{3.3}) over $t \in E_{n-1}$, by using (\ref{1.10}) and (\ref{3.1}), we obtain the required result.
\end{proof}

Now we discuss our main result which  are some generalized inequalities
involving $\beta_{\zeta,k}^{\left(a_{n},b_{n}\right)}{(\phi;\zeta)}$.

\begin{thm}
\label{thm21}
Let $ a_{n} >0$, $b_{n} - a_{n} > 0$, $\eta>0$ and $\zeta \geq 0$.
Then
\begin{equation}
\label{4.1}
\displaystyle
\left(\beta_{\zeta,k}^{\left(a_{n},b_{n}\right)}{(\phi+\psi;\eta)}\right)^{2} \leq
\beta_{\zeta,k}^{\left(a_{n},b_{n}\right)}{(2\phi;\eta)}\beta_{\zeta,k}^{\left(a_{n},b_{n}\right)}{(2\psi;\eta)}
\end{equation}
holds for any $\phi, \psi \in (0,\infty)^{n}$.
The real valued function $\beta_{\zeta,k}^{\left(a_{n},b_{n}\right)}{(\phi;\eta)}$ is convex on $(0,\infty)^{n}$.
\end{thm}

\begin{proof}
Let $\phi = \left(\phi_{1}, \ldots, \phi_{n}\right)$, $\psi= \left(\psi_{1}, \ldots, \psi_{n}\right)$ and
\begin{equation*}
\displaystyle
\lambda(t)= \ {}_{1}F_{1,k} \left(a_{n}, b_{n},-\frac {\eta^{k}}{k\pi(t)}-\zeta^{k}\frac {\pi(t)}{\eta^{k}} \right).
\end{equation*}
As $\lambda(t)\geq 0$,  we can write
\begin{equation*}
\displaystyle
\left(\beta_{\zeta,k}^{\left(a_{n},b_{n}\right)}{(\phi+\psi;\eta)}\right)^{2} =
\frac {1}{k^{n-1}}\left(\int\limits_{E_{n-1}}
\left(\prod_{i=1}^{n}t_{i}^{\frac {\phi_{i}}{k}-\frac {1}{2}}\right)
\left(\lambda(t)\right)^{\frac {1}{2}}
\left(\prod_{i=1}^{n}t_{i}^{\frac {\psi_{i}}{k}-\frac {1}{2}}
(\lambda(t))^{\frac {1}{2}}\right)dt\right)^{2}.
\end{equation*}
By  the Cauchy-Schwartz inequalities, (\ref{4.1}) is  obtained.

Now we have an interesting remark that (\ref{4.1}) is equivalent to
\begin{equation*}
\displaystyle
\beta_{\zeta,k}^{\left(a_{n},b_{n}\right)}{\left(\frac {\phi+\psi}{2};\eta\right)}
\leq \left(\beta_{\zeta,k}^{\left(a_{n},b_{n}\right)}{(\phi;\eta)}\beta_{\zeta,k}^{\left(a_{n},b_{n}\right)}{(\phi;\eta)}\right)^{\frac {1}{2}}.
\end{equation*}
By using arithmetic-geometric mean inequality $\sqrt{\phi \psi}\leq\frac {1}{2}\phi+\frac {1}{2}\psi$,
\begin{equation*}
\displaystyle
\beta_{\zeta,k}^{\left(a_{n},b_{n}\right)}{\left(\frac {\phi+\psi}{2};\eta\right)}\leq
\frac {1}{2}\beta_{\zeta,k}^{\left(a_{n},b_{n}\right)}{(\phi;\eta)}+\frac {1}{2}\beta_{\zeta,k}^{\left(a_{n},b_{n}\right)}{(\phi;\eta)}.
\end{equation*}
This expression shows that $\phi\longmapsto\beta_{\zeta,k}^{\left(a_{n},b_{n}\right)}$ is mid convex.
In addition of the fact that $\phi\longmapsto\beta_{\zeta,k}^{\left(a_{n},b_{n}\right)}$ has the  continuous property,
these facts  ensure $\phi\longmapsto\beta_{\zeta,k}^{\left(a_{n},b_{n}\right)}$ is convex.
The  proof is complete.
\end{proof}

\begin{lemma}
Let $\phi  > 0$, $\psi > 0$.
Then we have a real valued function $w \longmapsto  \ {}_{1}F_{1,k}  \left(a_{n};b_{n};w\right)$
which is differentiable on $\Re$ and
\begin{equation}
\label{4.2}
\displaystyle
\frac {d}{dw} \ {}_{1}F_{1,k} \left(a_{n};b_{n};w\right) =
\frac {a_{n}}{b_{n}} \ {}_{1}F_{1,k} \left(a_{n}+k;b_{n}+k;w\right).
\end{equation}
If $w = 0$, we have
\begin{equation}
\label{4.3}
\displaystyle
\frac {d}{dw} \ {}_{1}F_{1,k} \left(a_{n};b_{n};0\right) = \frac {a_{n}}{b_{n}}.
\end{equation}
\end{lemma}

Now we  state the following result.

\begin{thm}
\label{thm22}
Let $ a_{n} > 0$, $b_{n}-a_{n} >0$, $\eta >0$ and $\zeta \geq 0 $.
Then
\begin{equation*}
\displaystyle
\beta_{k}^{\left(a_{n},b_{n}\right)}{(\phi;\eta)} -
\frac {a_{n}}{b_{n}}\frac {\zeta^{k}}{\eta^{k}}\beta_{k}^{\left(a_{n}+1,b_{n}+1\right)}{(\phi+ke;\eta)}
\leq \beta_{\zeta,k}^{\left(a_{n},b_{n}\right)}{(\phi;\eta)} \leq
\beta_{k}^{\left(a_{n},b_{n}\right)}{(\phi;\eta)} \leq  \beta_{k}{(\phi)}
\end{equation*}
is valid for all $\phi \in (0,\infty)^{n}$, and $e=:\left(1,1, \ldots, 1\right)$.
\end{thm}

\begin{proof}
Using Remark \ref{11},  (\ref{3.1}) and (\ref{3.2}),
we  have $ \beta_{\zeta,k}^{\left(a_{n}, b_{n}\right)}{(\phi;\eta)} \leq \beta_{k}^{\left(a_{n},b_{n}\right)}{(\phi;\eta)} \leq \beta_{k}{(\phi)}$.
To prove the next part of the inequality
\begin{equation}
\label{4.5}
\displaystyle
\beta_{k}^{\left(a_{n},b_{n}\right)}{(\phi;\eta)} -
\frac {a_{n}}{b_{n}}\frac {\zeta^{k}}{\eta^{k}}\beta_{k}^{\left(a_{n}+1,b_{n}+1\right)}{(\phi+ke;\eta)}
\leq \beta_{\zeta,k}^{\left(a_{n},b_{n}\right)}{(\phi;\eta)}
\end{equation}
for fixed  $a_{n}, b_{n} - a_{n} > 0$ the map $b\longmapsto \ {}_{1}F_{1,k}\left(a_{n};b_{n};b\right)$ is convex on $\Re$.
Also we have $ f:\Re \longrightarrow \Re$ a convex function, differentiable at $c_{0}$,
so $ f(c)\geq f \left(c_{0}\right) + \left(c-c_{0}\right)\acute{f} \left(c_{0}\right)$.
Applying this and utilizing (\ref{4.2}), we have
\begin{equation}
\label{4.6}
\displaystyle
{}_{1}F_{1,k} \left(a_{n};b_{n};-\frac {\eta^{k}}{k\pi(t)}-\zeta^{k}\frac {\pi(t)}{\eta^{k}}\right)
\geq \ {}_{1}F_{1,k}{\left(a_{n};b_{n};-\frac {\eta^{k}}{k\pi(t)}\right)} -
\frac {a_{n}}{b_{n}}\frac {\zeta^{k}}{\eta^{k}}\pi(t)_{1}F_{1,k}\left(a_{n}+1;b_{n}+1;-\frac {\eta^{k}}{k\pi(t)}\right).
\end{equation}
Next multiply (\ref{4.6})  with $\prod_{i=1}^{n}t^{\frac {\phi_{i}}{k}-1}$ and
integrate over $t \in E_{n-1}$, using (\ref{3.1}) and (\ref{3.2}), we have (\ref{4.5}).
\end{proof}

\begin{proposition}
Let $ a_{n} > 0$, $b_{n} - a_{n} >0$, $\eta >0$ and $\zeta \geq 0$.
For $\phi \in (1,\infty)^{n}$, we have
\begin{equation*}
\displaystyle
\beta_{\zeta,k}^{\left(a_{n},b_{n}\right)}{(\phi;\eta)}\geq
\beta{(\phi,k)}-\frac {a_{n}}{b_{n}}\frac {\eta^{k}}{k} \beta_{k}{(\phi-ke)}-
\frac {a_{n}}{b_{n}}\frac {\zeta^{k}}{\eta^{k}}\beta_{k}{(\phi+ke)}
\end{equation*}
with $e$ defined as  earlier.
\end{proposition}

\begin{proof}
Using (\ref{4.3}), we have
\begin{equation}
\label{4.8}
\displaystyle
{}_{1}F_{1,k} \left(a_{n};b_{n};-\frac {\eta^{k}}{k\pi(t)}-\zeta^{k}\frac {\pi(t)}{\eta^{k}}\right)
\geq \ {}_{1}F_{1,k}{\left(a_{n};b_{n};0\right)} +
\frac {a_{n}}{b_{n}} \left(-\frac {\eta^{k}}{k\pi(t)}-\zeta^{k}\frac {\pi(t)}{\eta^{k}}\right).
\end{equation}
Multiplying (\ref{4.8}) with $\prod_{i=1}^{n}t^{\frac {\phi_{i}}{k}-1}$,
the rest of the proof is similar to that of Theorem \ref{thm22}.
\end{proof}

\begin{lemma}
\label{lem22}
We have
\begin{equation*}
\displaystyle
\sup_{t \in E_{n-1}} \pi (t) =x^{-x}.
\end{equation*}
\end{lemma}

\begin{proof}
Take $n$ positive real numbers $c_{1},c_{2}, \ldots, c_{n}$ and use
the arithmetic-geometric mean inequality     $\sqrt[n]{c_{1} \cdots  c_{n}} \leq \frac {c_{1}+ \cdots +c_{n}}{n}$.
Letting  $c_{1}=t_{1}, c_{2}=t_{2}, \ldots, c_{n-1}=t_{n-1}, c_{n}=t_{n}=:1-t_{1}- \cdots  -t_{n-1}$ in the  inequality,
we have     $\sqrt[n]{\pi(t)} \leq \frac {1}{n}$ or    $\pi(t)\leq n^{-n}$
for any $ t\in E_{n-1}$.
This inequality is an   equality for $\left(t_{1}, \ldots, t_{n-1}\right)=\left(\frac {1}{n}, \ldots, \frac {1}{n}\right)$.
\end{proof}

We now prove a refinement of the inequality  $\beta_{\zeta,k}^{\left(a_{n},b_{n}\right)}{(\phi;\eta)} \leq \beta_{k}{(\phi)}$.

\begin{thm}
\label{thm23}
Let $a_{n} > 0$, $b_{n}-a_{n} > 0$, $\eta >0$, $\zeta \geq 0$ and $\phi \in (0,\infty)^{n}$.
Also assume that $\eta\geq \sqrt{\zeta}$.
Then,
\begin{equation}
\label{4.11}
\displaystyle
\frac {\beta_{\zeta,k}^{\left(a_{n},b_{n}\right)}{(\phi;\eta)}}{\beta_{k}{(\phi)}}
\leq \ {}_{1}F_{1,k}\left(a_{n};b_{n};-\frac {\eta^{k}}{k}n^{n}-\frac {\zeta^{k}}{\eta^{k}n^{n}}\right)
\leq \ {}_{1}F_{1,k}\left(a_{n};b_{n};-2\sqrt{\zeta}\right) \leq 1,
\end{equation}
where $\beta_{k}{(\phi)}$ is as   defined in (\ref{1.10}).
\end{thm}

\begin{proof}
First we consider the inequality
\begin{equation*}
\displaystyle
{}_{1}F_{1,k}\left(a_{n};b_{n};-\frac {\eta^{k}}{k}n^{n}-\frac {\zeta^{k}}{\eta^{k}n^{n}}\right)
\leq \ {}_{1}F_{1,k}\left(a_{n};b_{n};-2\sqrt{\zeta}\right) \leq 1.
\end{equation*}
We are able to verify it easily by considering
\begin{equation*}
\displaystyle
-\frac {\eta^{k}}{k}n^{n}-\frac {\zeta^{k}}{\eta^{k}n^{n}} \leq -2\sqrt{\zeta} \leq 0,
\quad
\displaystyle
\ {}_{1}F_{1,k}\left(a_{n};b_{n};0\right)=1
\end{equation*}
and $ w \longmapsto \ {}_{1}F_{1,k}\left(a_{n};b_{n};w\right)$ is a   real valued function and also increasing.
We will prove the next part of the inequality
\begin{equation}
\label{4.12}
\displaystyle
\frac {\beta_{\zeta,k}^{\left(a_{n},b_{n}\right)}{(\phi;\eta)}}{\beta_{k}{(\phi)}} \leq
\ {}_{1}F_{1,k}\left(a_{n};b_{n};-\frac {\eta^{k}}{k}n^{n}-\frac {\zeta^{k}}{\eta^{k}n^{n}}\right).
\end{equation}
By using (\ref{3.1}), we obtain
\begin{equation}
\label{4.13}
\displaystyle
\beta_{\zeta,k}^{\left(a_{n},b_{n}\right)}{(\phi;\eta)} \leq
\beta_{k}{(\phi)}
\sup_{E_{n-1}} \ {}_{1}F_{1,k} \left(a_{n};b_{n};-\frac {\eta^{k}}{k\pi(t)}-\zeta^{k}\frac {\pi(t)}{\eta^{k}}\right).
\end{equation}
The  supremum exists and
\begin{equation}
\label{4.14}
\displaystyle
{}_{1}F_{1,k}\left(a_{n};b_{n};-\frac {\eta^{k}}{k\pi(t)}-\zeta^{k}\frac {\pi(t)}{\eta^{k}}\right) =
\frac {\Gamma \left(b_{n}\right)}{k\Gamma\left(a_{n}\right)\Gamma\left(b_{n}-a_{n}\right)}
\int\limits_{0}^{1}b^{\frac {a_{n}}{k}-1}
(1-b)^{\frac {b_{n}-a_{n}}{k}-1}\exp \left[ b\left(-\frac {\eta^{k}}{k\pi(t)}-\zeta^{k}\frac {\pi(t)}{\eta^{k}}\right) \right] db.
\end{equation}
It sufficient to find an upper bond of
$t\mapsto-\frac {\eta^{k}}{k\pi(t)}-\zeta^{k}\frac {\pi(t)}{\eta^{k}}$.
Since $b \in [0,1]$ and the exponential function is increasing,
if $\eta \geq \sqrt{\zeta}$ and $0 < u < 1$, then
\begin{equation}
\label{4.15}
\displaystyle
\sup_{0\leq s\leq u}\left(-\frac {\eta^{k}}{ks}-\zeta^{k}\frac {s}{\eta^{k}}\right)= -\frac {\eta^{k}}{ku}-\zeta^{k}\frac {u}{\eta^{k}}.
\end{equation}
By  Lemma \ref{lem22}, $\pi(t) \leq x^{-x}$ and so
$E_{n-1}\subset \left\{ \left( t_{1}, \ldots, t_{n-1} \right): 0 \leq \pi(t) \leq x^{-x} \right\}=:G$.
This with (\ref{4.15}) implies that
\begin{equation*}
\displaystyle
\sup_{t\epsilon E_{n-1}}
\left(-\frac {\eta^{k}}{k\pi(t)}-\zeta^{k}\frac {\pi(t)}{\eta^{k}}\right)
\leq \sup_{t \in G}
\left(-\frac {\eta^{k}}{k\pi(t)}-\zeta^{k}\frac {\pi(t)}{\eta^{k}} \right) =
\sup_{0 \leq s \leq x^{-x}}  \left(-\frac {\eta^{k}}{ks}-\zeta^{k}\frac {s}{\eta^{k}}\right) =
{\frac {-\eta^{k}}{k}}n^{n}-\frac {\zeta^{k}}{\eta^{k}n^{n}}.
\end{equation*}
Substituting into  (\ref{4.14}) and then using  (\ref{4.13}), we have (\ref{4.12}).
\end{proof}

\begin{corollary}
Under the assumptions of  Theorem \ref{thm23} for any $ z > 0$, we have
\begin{equation}
\label{4.16}
\displaystyle
\int\limits_{0}^{\infty}\eta^{z-1}\beta_{k}^{\left(a_{n},b_{n}\right)}{(\phi;\eta)}d\eta \leq
\beta_{k} {(\phi)}n^{\left(-n^{-n}\right)^{\frac {z}{k}}} \Gamma_{k}^{\left(a_{n},b_{n}\right)}{(z)},
\end{equation}
where $\Gamma_{k}^{\left(a_{n},b_{n}\right)}{(z)}$ is defined in (\ref{2.4}).
If $z = 1$
\begin{equation*}
\displaystyle
\int\limits_{0}^{\infty}\beta_{k}^{\left(a_{n},b_{n}\right)}{(\phi;\eta)}d\eta \leq
\beta_{k} {(l)}{n^{-n}}^{\frac {1}{k}} \Gamma^{\left(a_{n},b_{n}\right)}{(1)}.
\end{equation*}
\end{corollary}

\begin{proof}
Taking $\zeta=0$ in (\ref{4.11}), we have
\begin{equation*}
\displaystyle
\beta_{k}^{\left(a_{n},b_{n}\right)}{(\phi;\eta)} \leq \beta_{k}{(\phi)}
\ {}_{1}F_{1,k}\left(a_{n};b_{n};-\frac {\eta^{k}}{k}n^{n}\right).
\end{equation*}
Integrating over $\eta \in (0,\infty)$ and multiplying
by   $\eta^{z-1}$, we obtain (\ref{4.16}) after a  simple change of variables.
\end{proof}

Next we state a  result which deals with the lower bound of $\beta_{\zeta,k}^{\left(a_{n},b_{n}\right)}{(\phi;\eta)}$.

\begin{thm}
\label{thm24}
Let $ a_{n} > 0$, $b_{n} - a_{n} > 0$, $\eta > 0$, $\zeta \geq 0$ and $\phi \in (0,\infty)$.
Then,
\begin{equation}
\label{4.18}
\displaystyle
\frac {\beta_{\zeta,k}^{\left(a_{n},b_{n}\right)}{(\phi;\eta)}}{\beta_{k}{(\phi;\eta)}} \geq
\ {}_{1}F_{1,k}\left(b_{n}-a_{n};b;\frac {\eta^{k}n^{n}}{k}\right)
\exp\left(-\frac {\zeta^{k}}{\eta^{k}n^{n}}\right),
\end{equation}
where $ \beta_{k}{(\phi;\eta)}$ is  defined in (\ref{1.14}).
\end{thm}

\begin{proof}
Using (\ref{1.14}), we have
\begin{equation*}
\displaystyle
\beta_{\zeta,k}^{\left(a_{n},b_{n}\right)}{(\phi;\eta)}=\frac {1}{k^{n-1}}\int\limits_{E_{n-1}}
\prod_{i=1}^{n}t_{i}^{\frac {\phi_{i}}{k}-1}
\exp \left(-\frac {\eta^{k}}{k\pi(t)}-\zeta^{k}\frac {\pi(t)}{\eta^{k}}\right)
\ {}_{1}F_{1,k}\left(b_{n}-a_{n};b_{n};\frac {\eta^{k}}{k\pi(t)}+\zeta^{k}\frac {\pi(t)}{\eta^{k}}\right)dt,
\end{equation*}
which can be  rewritten as
\begin{equation}
\label{4.19}
\displaystyle
\beta_{\zeta,k}^{\left(a_{n},b_{n}\right)}{(\phi;\eta)} =
\frac {1}{k^{n-1}}\int\limits_{E_{n-1}}
\left(\prod_{i=1}^{n}t_{i}^{\frac {\phi_{i}}{k}-1} e^{-\frac {\eta^{k}}{k\pi(t)}}\right)
e^{-\frac {\zeta^{k}\pi(t)}{\eta^{k}}}
\ {}_{1}F_{1,k}\left(b_{n}-a_{n};b_{n};\frac {\eta^{k}}{k\pi(t)}+\zeta^{k}\frac {\pi(t)}{\eta^{k}}\right)dt.
\end{equation}
Using Lemma \ref{lem22}, we have $e^{-\frac {\zeta^{k}\pi(t)}{\eta^{k}}} \geq e^{-\frac {\zeta^{k}}{\eta^{k}n^{n}}}$
and $\frac {\eta^{k}}{k\pi(t)} + \zeta^{k}\frac {\pi(t)}{\eta^{k}} \geq \frac {\eta^{k}}{k}n^{n}$.
Using this in \eqref{4.19} and with the help of \eqref{1.14}, we have \eqref{4.18}.
\end{proof}

\begin{corollary}
Under the assumptions of  Theorem \ref{thm24}, for any $z > 0$,  we have
\begin{equation}
\label{4.20}
\displaystyle
\int\limits_{0}^{\infty} r^{z-1}\beta_{\zeta,k}^{\left(a_{n},b_{n}\right)}{(\phi;\eta)}dr \geq
\eta^{z} \left(\frac {n^{n}}{k}\right)^{\frac {z}{k}}\Gamma_{k}{(z)}\beta_{k}{(\phi;\eta)}
\ {}_{1}F_{1,k}\left(b_{n}-a_{n};b_{n};\frac {\eta^{k}{k}}n^{n}\right).
\end{equation}
If $ z = 1$ then
\begin{equation}
\label{4.21}
\displaystyle
\int\limits_{0}^{\infty}
\beta_{\zeta,k}^{\left(a_{n},b_{n}\right)}{(\phi;\eta)}dr \geq
\eta^{z} \left(\frac {n^{n}}{k}\right)^{\frac {z}{k}} \beta_{k}{(\phi;\eta)}\Gamma_{k}{(1)}
\ {}_{1}F_{1,k} \left(a_{n}-b_{n};b_{n};\frac {\eta^{k}{k}}n^{n}\right).
\end{equation}
\end{corollary}

\begin{proof}
Multiplying (\ref{4.18})  with $r^{z-1}$ and integrating over $r \in (0,\infty)$, we have
\begin{equation*}
\displaystyle
\int\limits_{0}^{\infty} r^{z-1}{\beta_{\zeta,k}^{\left(a_{n},b_{n}\right)}{(\phi;\eta)}}dr
\geq \beta_{k}{(\phi;\eta)} \ {}_{1}F_{1,k}\left(a_{n}-b_{n};b_{n};\frac {\eta^{k}}{k}n^{n}\right)
\int\limits_{0}^{\infty}r^{z-1} \exp \left(  -\frac {\zeta^{k}}{\eta^{k}n^{n}} \right)dr.
\end{equation*}
Setting  $t = \left(\frac {kr^{k}}{\eta^{k}n^{n}}\right)^{\frac {1}{k}}$ and making simplification, we obtain (\ref{4.20}) and then (\ref{4.21}).
\end{proof}

Our next result is as follows.

\begin{thm}
\label{thm25}
Let $ a_{n} > 0$, $b_{n} - a_{n} > 0$, $\eta > 0$, $\zeta \geq 0$ and $\phi \in (0,\infty)$.
Then,
\begin{equation}
\label{4.23}
\displaystyle
\frac {\beta_{\zeta,k}^{\left(a_{n},b_{n}\right)}{(\phi;\eta)}}{\beta_{k}{(\phi;\eta)}}
\geq \ {}_{1}F_{1,k}\left(a_{n};b_{n};-\frac {\zeta^{k}}{\eta^{k}n^{n}}\right).
\end{equation}
\end{thm}

\begin{proof}
Using (\ref{1.8}), we have
{\tiny\begin{equation*}
\displaystyle
{}_{1}F_{1,k}\left(a_{n};b_{n};-\frac {\eta^{k}}{k\pi(t)}-\zeta^{k}\frac {\pi(t)}{\eta^{k}}\right) =
\frac {\Gamma \left(b_{n}\right)}{k\Gamma \left(a_{n}\right)\Gamma \left(b_{n}-a_{n}\right)}
\int\limits_{0}^{1}b^{\frac {a_{n}}{k}-1}(1-b)^{\frac {b_{n}-a_{n}}{k}-1}
\exp \left(-b\frac {\eta^{k}}{k\pi(t)}\right)
\exp\left(-b\zeta^{k}\frac {\pi(t)}{\eta^{k}}\right)db.
\end{equation*}}
Using Lemma \ref{lem22},  $\exp \left(-b\frac {\eta^{k}}{k\pi(t)}\right) \geq \exp \left(-\frac {\eta^{k}}{k\pi(t)}\right)$
and $\exp \left(-b\zeta^{k}\frac {\pi(t)}{\eta^{k}}\right) \geq \exp \left(-b\frac {\zeta^{k}}{\eta^{k}n^{n}}\right)$.
Applying these in
\begin{equation*}
\displaystyle
\beta_{\zeta,k}^{\left(a_{n},b_{n}\right)}{(\phi;\eta)} =
\frac {1}{k^{n-1}} \int\limits_{E{_{n-1}}}\prod_{i=1}^{n}t_{i}^{\frac {\phi_{i}}{k}-1}
\ {}_{1}F_{1,k}\left(a_{n};b_{n};-\frac {\eta^{k}}{k\pi(t)}-\zeta^{k}\frac {\pi(t)}{\eta^{k}} \right)dt,
\end{equation*}
and   using the uniform convergence of the involved integrals, we have
{\tiny\begin{equation*}
\displaystyle
\beta_{\zeta,k}^{\left(a_{n},b_{n}\right)}{(\phi;\eta)} \geq\frac {1}{k^{n-1}}
\frac {\Gamma(b_{n})}{k\Gamma \left(a_{n}\right)\Gamma\left(b_{n}-a_{n}\right)}
\int\limits_{0}^{1}b^{\frac {a_{n}}{k}-1}(1-b)^{\frac {b_{n}-a_{n}}{k}-1}
\exp \left(-b\frac {\zeta^{k}}{\eta^{k}n^{n}}\right)db
\int\limits_{E_{n-1}}\prod_{i=1}^{n}t_{i}^{\frac {\phi_{i}}{k}-1} \exp\left(-\frac {\eta^{k}}{k\pi(t)}\right)dt.
\end{equation*}}
Hence, \eqref{4.23}.
\end{proof}

\begin{corollary}
Under the assumptions of  Theorem \ref{thm24}, for any $z > 0$,  we have
\begin{equation*}
\displaystyle
\int\limits_{0}^{\infty}
\zeta^{z-1} \beta_{\zeta,k}^{\left(a_{n},b_{n}\right)}{(\phi;\eta)}d\zeta \geq
\eta^{z}\left(\frac {n^{n}}{k}\right)^\frac {z}{k}\beta_{k}{(\phi;\eta)}\Gamma_{k}^{\left(a_{n},b_{n}\right)}{(z)},
\end{equation*}
where $\Gamma_{k}^{(c,d)}{(z)}$ is defined  in \eqref{2.4}.
If $ z = 1$ then
\begin{equation*}
\displaystyle
\int\limits_{0}^{\infty} \beta_{\zeta,k}^{\left(a_{n},b_{n}\right)}d\zeta \geq
\eta \left(\frac {n^{n}}{k}\right)^\frac {1}{k} \beta_{k}{(\phi;\eta)}\Gamma_{k}^{\left(a_{n},b_{n}\right)}{(1)}.
\end{equation*}
\end{corollary}

\begin{proof}
Multiplying \eqref{4.23} with $v^{z-1}$ and integrating over $v \in (0,\infty)$, we obtain
\begin{equation*}
\displaystyle
\int\limits_{0}^{\infty}
\zeta^{z-1} \beta_{\zeta,k}^{\left(a_{n},b_{n}\right)}{(\phi;\eta)}d\zeta \geq
\beta_{k}{(\phi;\eta)}\int\limits_{0}^{\infty}\zeta^{z-1}
\ {}_{1}F_{1,k}\left(a_{n};b_{n};-\frac {\zeta^{k}}{\eta^{k}n^{n}}\right) dv.
\end{equation*}
Setting  $ \zeta = \eta t\frac {\left(n^{n}\right)^{\frac {1}{k}}}{k^{\frac {1}{K}}}$ and simplifying yields the desired result.
\end{proof}

\subsection{Generalized beta $k$  function of the second kind}

In this section, we discuss the generalized beta $k$ function of the second kind.
Before describing its representation, we give some notation  used throughout the section.
Let   $\phi = \left(\phi_{1}, \ldots, \phi_{n}\right)$,
$\eta = \left(\eta_{1}, \ldots, \eta_{n} \right)$,
$p = \left(p_{1}, \ldots, p_{n}\right)$,
$q = \left(q_{1}, \ldots, q_{n} \right)$ and  $\zeta = \left(\zeta_{1}, \ldots, \zeta_{n}\right)$.

\begin{definition}
Let $\phi, p, q, \eta \in (0,\infty)^{n}$ and $  \zeta\in [0,\infty)^{n}$, we
define a   generalized beta $k$ function as
\begin{equation}
\label{5.1}
\displaystyle
\beta_{\zeta,k}^{(p,q)}{(\phi;\eta)}=:\frac {1}{k^{n-1}}\int\limits_{{E_{n-1}}}
\prod_{i=1}^{n}t_{i}^{\frac {\phi_{i}}{k}-1}
\ {}_{1}F_{1,k}\left(p_{i};q_{i};-\frac {\eta_{i}^{k}}
{kt_{i}}-\zeta_{i}^{k}\frac {t_{i}}{f_{i}^{k}}\right)dt,
\end{equation}
where $dt = \left(dt_{1}, \ldots, dt_{n-1}\right)$ and $t_{n} = 1-\sum_{i=1}^{n-1}t_{i}$.
If $ \zeta = 0$ then
\begin{equation}
\label{5.2}
\displaystyle
\beta^{(p,q)}{(\phi;\eta)}=:\frac {1}{k^{n-1}}
\int\limits_{{E_{n-1}}} \prod_{i=1}^{n}t_{i}^{\frac {\phi_{i}}{k}-1}
\ {}_{1}F_{1,k} \left(p_{i};q_{i};-\frac {{\eta_{i}}^{k}}{kt_{i}}\right) dt.
\end{equation}
If $p=q$ and $n=2$ then (\ref{5.2}) becomes similar to  (\ref{1.4}).
If $p = q$ and $\zeta =0$ then (\ref{5.1}) is exactly (\ref{1.15}).
\end{definition}

\begin{remark}
(\ref{5.1}) is well defined because following the
uniform convergence of the stated series in (\ref{1.7}) we are able to
interchange series and integral defined in (\ref{5.1}).
Further such integrals are uniformly convergent in any compact set included in the interior of $E_{n-1}$.
This allows for  differentiation and limit under the
integral sign of (\ref{5.1}).
We may state  $\lim_{\zeta \rightarrow 0}\beta_{\zeta,k}^{(p,q)}{(\phi;\eta)} = \beta_{k}^{(p,q)}{(\phi;\eta)}$.
\end{remark}

\begin{proposition}
For any $p, \phi, \eta \in (0,\infty)^{n}$, we have
\begin{equation*}
\displaystyle
\beta_{k}{(\phi;\eta)}= \lim_{\zeta \rightarrow 0}\beta_{\zeta,k}^{(p;p)}{(\phi;\eta)}=:\beta_{k}^{(p,p)}{(\phi;\eta)}.
\end{equation*}
\end{proposition}

\begin{proof}
Using (\ref{1.7}) and (\ref{5.2}), we have for  $i= 1,2, \ldots, n$
\begin{equation*}
\displaystyle
{}_{1}F_{1,k}\left(p_{i};p_{i};-\frac {\eta_{i}^{k}}{kt_{i}}\right)=
\sum_{i=1}^{\infty}\frac {\left(-\frac {\eta_{i}^{k}}{kt_{i}}\right)^{m}}{m!}    = e^{-\frac {\eta_{i}^{k}}{kt_{i}}}.
\end{equation*}
Combining  (\ref{1.15}) and (\ref{5.1}) gives the required result.
\end{proof}

Next we  mention some inequalities involving  the described beta function.
Our first result is   related to  convexity of $\beta_{\zeta,k}^{(p,q)}{(\phi;\eta)}$.

\begin{thm}
Let $ \eta, p, q - p \in (0, \infty) $ and $ \zeta \in [0,\infty)$.
Then,
\begin{equation*}
\displaystyle
\left(\beta_{\zeta,k}^{(p,q)}{(\phi + \psi;\eta)}\right)^{2} \leq
\beta_{\zeta,k}^{(p,q)}{(2\phi;\eta)}\beta_{\zeta,k}^{(p,q)}{(2\psi;\eta)}
\end{equation*}
holds for any $\phi, \psi \in (0,\infty)^{n}$.
Therefore, $\beta_{\zeta,k}^{(p,q)}{(\phi;\eta)}$ is convex on $(0,\infty)^{n}$.
\end{thm}

\begin{proof}
The proof is similar to that of Theorem \ref{thm21}.
\end{proof}

\begin{thm}
We have
\begin{equation*}
\displaystyle
\frac {\beta_{\zeta,k}^{(p,q)}{(\phi;\eta)}}{\beta_{k}(l)}
\leq \prod_{i=1}^{n} e^{-\frac {\eta_{i}^{k}}{k}}
\ {}_{1}F_{1,k}\left(q_{i}-p_{i};q_{i};\frac {\eta_{i}^{k}}{k}+\frac {\zeta_{i}^{k}}{\eta_{i}^{k}}\right),
\end{equation*}
where $\beta_{k}{(\phi)}$ is  defined in (\ref{1.10}).
\end{thm}

\begin{proof}
Since the real valued function $l \mapsto\  {}_{1}F_{1,k}\left(c_{n};d_{n};l\right)$ is
increasing and $0 < t_{i} \leq 1$, $ i = 1, \ldots, n$, we have
\begin{equation*}
\displaystyle
\beta_{\zeta,k}^{(p,q)}{(\phi;\eta)} \leq
\frac {1}{k^{n-1}}\int\limits_{{E_{n-1}}}\prod_{i=1}^{n}t_{i}^{\frac {\phi_{i}}{k}-1}
\ {}_{1}F_{1,k}\left(p_{i};q_{i};-\frac {\eta_{i}^{k}}{k}-\zeta_{i}^{k}\frac {t_{i}}{\eta_{i}^{k}}\right)dt.
\end{equation*}
Using (\ref{1.9}),
\begin{equation}
\label{6.2}
\displaystyle
\beta_{\zeta,k}^{(p,q)}{(\phi;\eta)} \leq \frac {1}{k^{n-1}}
\int\limits_{{E_{n-1}}}\prod_{i=1}^{n}t_{i}^{\frac {\phi_{i}}{k}-1}
\exp \left({-\frac {\eta_{i}^{k}}{kt_{i}}}-\frac {\zeta^{k}t_{i}}{\eta_{i}^{k}}\right)
\ {}_{1}F_{1,k}\left(q_{i}-p_{i};q_{i};\frac {\eta_{i}^{k}}{k}+\zeta_{i}^{k}\frac {t_{i}}{\eta_{i}}^{k}\right)dt.
\end{equation}
It is clear that
$ \exp\left(\frac {-\eta_{i}^{k}}{k}\right)-\zeta_{i}^{k}\frac {t_{i}}{\eta_{i}^{k}} \leq e^{\frac {-\eta_{i}^{k}}{k}}$
and $\frac {\eta_{i}^{k}}{k}+\zeta_{i}^{k}\frac {t_{i}}{\eta_{i}^{k}} \leq \frac {\eta_{i}^{k}}{k}+\frac {\zeta_{i}^{k}}{\eta_{i}^{k}}$.
Using these in (\ref{6.2}),  we obtain the desired result.
\end{proof}

Now we provide a  lower bond of the described function.

\begin{thm}
We have
\begin{equation}
\label{6.7}
\displaystyle
\frac {\beta_{\zeta,k}^{(p,q)}{(\phi;\eta)}}{\beta_{k}(\phi;\eta)}\geq
\prod_{i=1}^{n}e^{-\frac {\zeta_{i}^{k}}{\eta_{i}^{k}}}
\ {}_{1}F_{1,k}\left(q_{i}-p_{i};q_{i};m_{i}\right),
\end{equation}
where $\beta{(\phi;\eta)} $ is defined in (\ref{1.15}) and
$m_{i}=\max \left(2\sqrt{\zeta_{i}^{k}},\eta_{i}^{k}\right)$.
If $ \eta_{i}^{k} \geq \sqrt{\zeta_{i}^{k}}$, $i=1, \ldots, n$, then (\ref{6.7}) can   be refined as
\begin{equation}
\label{6.8}
\displaystyle
\frac {\beta_{\zeta,k}^{(p,q)}{(\phi;\eta)}}{\beta_{k}(\phi;\eta)} \geq
\prod_{i=1}^{n}e^{-\frac {\zeta_{i}^{k}}{\eta_{i}^{k}}}
\ {}_{1}F_{1,k}\left(q_{i}-p_{i};q_{i};\frac {\eta_{i}^{k}}{k}+\frac {\zeta_{i}^{k}}{\eta_{i}^{k}}\right).
\end{equation}
\end{thm}

\begin{proof}
Using (\ref{1.9}) and (\ref{5.1}), we obtain
\begin{equation}
\label{6.9}
\displaystyle
\beta_{\zeta,k}^{(p,q)}{(\phi;\eta)} =\frac {1}{k^{n-1}}
\int\limits_{E_{n-1}}
\left(\prod_{i=1}^{n}t_{i}^{\frac {\eta_{i}^{k}}{k}-1}
e^{-\frac {\eta_{i}^{k}}{kt_{i}}}\right)
e^{-\zeta_{i}^{k}\frac {t_{i}}{\eta_{i}^{k}}}
\ {}_{1}F_{1,k} \left(q_{i}-p_{i};q_{i}; \frac {\eta_{i}^{k}}{kt_{i}}+\zeta_{i}^{k}\frac {t_{i}}{\eta_{i}^{k}}\right) dt.
\end{equation}
It is easy to see that
$ \frac {\eta_{i}^{k}}{kt_{i}}+ \zeta_{i}^{k}\frac {t_{i}}{\eta_{i}^{k}} \geq 2\sqrt{\zeta_{i}^{k}}$,
also, further $0 < t_{i}<1$, and $h \mapsto \ {}_{1}F_{1}\left(a_{n};b_{n};c\right)$ is increasing.
We have $ e^{-\zeta_{i}^{k}\frac {t_{i}}{\eta_{i}^{k}}} \geq e^{-\frac {\zeta_{i}^{k}}{\eta_{i}^{k}}}$ and
${}_{1}F_{1,k} \left(q_{i}-p_{i};q_{i};\frac {\eta_{i}^{k}}{kt_{i}}+\zeta_{i}^{k}\frac {t_{i}}{\eta_{i}^{k}}\right)
\geq \ {}_{1}F_{1,k}\left(q_{i}-p_{i};q_{i};\max \left(2\sqrt{\zeta_{i}^{k}},\eta_{i}^{k}\right) \right)$.
Substituting into (\ref{6.9}) and using (\ref{1.15}), we have (\ref{6.7}).
$(\eta, \zeta)$ with $\eta \geq \sqrt{\zeta^{k}}$ implies $\inf
\left( \frac {\eta^{k}}{t}+\zeta^{k}\frac {t}{\eta^{k}}  \right) = \eta^{k}+\frac {\zeta^{k}}{\eta^{k}}$,
so   (\ref{6.8}) refines to   (\ref{6.7}).
The  proof is complete.
\end{proof}

We have the next result.

\begin{thm}
Let $ \phi, \eta, p, s-p \in (0,\infty)^{n}$.
For any $z =: \left(z_{1}, \ldots, z_{n} \right)$, we have
\begin{equation}
\label{6.10}
\displaystyle
\int\limits_{(0,\infty)^{n}}\prod_{i=1}^{n}\zeta_{i}^{z_{i}-1}
\beta_{\zeta,k}^{(p,q)}{(\phi;\eta)}dg \geq \beta_{k}{(\phi,\eta)}\Gamma_{k}{(z)}
\left(\prod_{i=1}^{n}\eta_{i}^{z_{i}}\frac {1}{k^{\frac {z_{i}}{k}}}\right)
\prod_{i=1}^{n}\ {}_{1}F_{1,k}\left(q_{i}-p_{i};q_{i};f_{i}^{k}\right),
\end{equation}
where $ d\zeta=: d\zeta_{1} \cdots  d\zeta_{n} $.
If   $z = e =(1, \ldots, 1)$ then
\begin{equation}
\label{6.11}
\displaystyle
\int\limits_{(0,\infty)^{n}}
\beta_{\zeta,k}^{(p,q)}{(\phi;\eta)}dg \geq \beta_{k}{(\phi,\eta)}
\left(\prod_{i=1}^{n}\eta_{i}^{z_{i}}\frac {1}{k^{\frac {z_{i}}{k}}}\right)
\prod_{i=1}^{n}\ {}_{1}F_ {1,k}\left(q_{i}-p_{i};q_{i};f_{i}^{k}\right).
\end{equation}
\end{thm}

\begin{proof}
Multiply (\ref{6.7}) with $\prod_{i=1}^{n}\zeta_{i}^{z_{i}-1}$  and
integrate over $ \zeta \in (0,  \infty)^{n} $ to obtain
\begin{equation*}
\displaystyle
\int\limits_{(0,\infty)^{n}}
\prod_{i=1}^{n}\zeta_{i}^{z_{i}-1}{\beta_{\zeta,k}^{(p,q)}{(\phi;\eta)}}dg \geq
\beta_{k}{(\phi;\eta)} \  {}_{1}F_{1,k}\left(q_{i}-p_{i};q_{i};\eta_{i}^{k}\right)
\int\limits_{(0,\infty)^{n}}\prod_{i=1}^{n}\zeta_{i}^{z_{i}-1}e^{-\frac {\zeta_{i}^{k}}{\eta_{i}^{k}}}dg.
\end{equation*}
Setting    $t = \left(\frac {k\zeta_{i}^{k}}{\eta_{i}^{k}}\right)^{\frac {1}{k}}, i=1, \ldots, n$, we have
\begin{equation*}
\displaystyle
\int\limits_{(0,\infty)^{n}}\prod_{i=1}^{n}\zeta_{i}^{z_{i}-1}{\beta_{\zeta,k}^{(p,q)}{(\phi;\eta)}}dg \geq \beta_{k}{(\phi;\eta)}
\prod_{i=1}^{n}{\eta_{i}^{z}\frac {1}{k^{\frac {z}{k}}}}
\ {}_{1}F_{1,k}\left(q_{i}-p_{i};q_{i};\eta_{i}^{k}\right)
\int\limits_{(0,\infty)^{n}}
\left(t_{i}\right)^{z_{i}-1}e^{-\frac {t_{i}^{k}}{k}}dt.
\end{equation*}
Hence (\ref{6.10}).
Taking $z_{i} = 1$, $i=1, \ldots, n$ in (\ref{6.10}), we have (\ref{6.11}).
The proof is complete.
\end{proof}

The following result may also be stated.

\begin{thm}
\label{thm210}
Let $ \phi, \eta, p, q-p \in (0,\infty)^{n}$ and $\zeta \in [0,\infty)$.
Then,
\begin{equation}
\label{6.14}
\displaystyle
\frac {\beta_{\zeta,k}^{(p,q)}{(\phi;\eta)}}{\beta_{k}(\phi;\eta)} \geq
\prod_{i=1}^{n}
\ {}_{1}F_{1,k}\left(p_{i};q_{i};-\frac {\zeta_{i}^{k}}{\eta_{i}^{k}}\right).
\end{equation}
\end{thm}

\begin{proof}
Using (\ref{1.8}), we have
\begin{equation}
\label{6.15}
\displaystyle
{}_{1}F_{1,k}\left(p_{i};q_{i};-\frac {\eta_{i}^{k}}{kt_{i}}-\zeta_{i}^{k}\frac {t_{i}}{\eta_{i}^{k}}\right) =
\frac {\Gamma\left(q_{i}\right)}{k\Gamma\left(p_{i}\right)-\Gamma\left(q_{i}-p_{i}\right)}
\int\limits_{0}^{1}u^{\frac {\eta_{i}}{k}-1}(1-u)^{\left(\frac {q_{i}-p_{i}}{k}-1\right)}
e^{-u\frac {\eta_{i}^{k}}{kt_{i}}}  e^{-u\zeta_{i}^{k}\frac {t_{i}}{\eta_{i}^{k}}}du.
\end{equation}
Since  $0 < t_{i} \leq 1$, $i=1, \ldots, n$ and $u \in [0,1]$, we have
$e^{-u\frac {\eta_{i}^{k}}{kt_{i}}} \geq e^{-\frac {\eta_{i}^{k}}{kt_{i}}}$
and   $e^{-\zeta_{i}^{k}u\frac {t_{i}}{\eta_{i}^{k}}} \geq e^{-u\frac {\zeta_{i}^{k}}{\eta_{i}^{k}}}$.
Using this in (\ref{6.15}) and then in (\ref{5.1}), we have
\begin{equation*}
\displaystyle
\beta_{\zeta,k}^{(p,q)}{(\phi;\eta)}\geq \frac {1}{k^{n-1}}\int\limits_{E_{n-1}}
\prod_{i=1}^{n}t_{i}^{\frac {\eta{_{i}}}{k}-1}
e^{-\frac {\eta_{i}^{k}}{k  t_{i}}}dt
\prod_{i=1}^{n}\frac {\Gamma \left(q_{i}\right)}{k\Gamma\left(p_{i}\right)-\Gamma\left(q_{i}-p_{i}\right)}
\int\limits_{0}^{1}u^{\frac {r_{i}}{k}-1}    (1-u)^{\frac {q_{i}-p_{i}}{k}-1}
e^{-u\frac {\zeta_{i}^{k}}{\eta_{i}^{k}}}du.
\end{equation*}
This with (\ref{1.15}) and (\ref{1.8})  yields (\ref{6.14}).
\end{proof}

Now  to our last result.

\begin{thm}
Let $ \phi, \eta, q, p-q \in (0,\infty)^{n}$.
For any $z=: \left(z_{1}, \ldots, z_{n}\right) \in (0,\infty)^{n}$, we have
\begin{equation*}
\displaystyle
\int\limits_{(0,\infty)^{n}}
\prod _{i=1}^{n}\zeta_{i}^{z_{{i}-1}}
\beta_{\zeta,k}^{(p,q)}{(\phi;\eta)}dg \geq
\beta_{k}{(\phi;\eta)}\Gamma_{k}^{(p,q)}{(z)}
\prod_{i=1}^{n}\left(\frac {\eta_{i}^{z_{i}}}{k^{\frac {z_{i}}{k}}}\right).
\end{equation*}
If $z = e =: (1, \ldots, 1)$ then
\begin{equation*}
\displaystyle
\int\limits_{(0,\infty)^{n}}
\beta_{\zeta,k}^{(p,q)}{(\phi;\eta)}dg \geq \beta_{k}{(\phi;\eta)}\Gamma_{k}^{(p,q)}{(e)}\prod_{i=1}^{n}\eta_{i}.
\end{equation*}
\end{thm}

\begin{proof}
Similar to the proof of Theorem \ref{thm210}.
Using (\ref{6.14}), we obtain
\begin{equation*}
\displaystyle
\beta_{\zeta,k}^{(p,q)}{(\phi;\eta)}\geq \beta_{k}{(\phi;\eta)}
\prod_{i=1}^{n}
\  {}_{1}F_{1,k} \left(p_{i};q_{i}; -\frac {\zeta_{i}^{k}}{\eta_{i}^{k}}\right).
\end{equation*}
Multiplying  with $\zeta_{i}^{z_{i}-1}$ and integrating over  $g \in (0, \infty)$,
we obtain our required result.
\end{proof}

\section{Conclusions}

In this note, we have presented generalized inequalities involving beta and gamma functions and their generalizations.
The note included basic representations of beta and gamma functions.
Some representations of beta and gamma functions involving confluent hypergeometric functions have been studied.
Some basic relations between  gamma and beta functions have been  provided.
Some refined inequalities involving extended beta function have been given.
We have also discussed upper and lower bounds for an extended beta function.

\subsection*{Ethics approval}

Not applicable.

\subsection*{Funding}

Not applicable.

\subsection*{Conflict of interest}

All of the   authors have no conflicts of interest.

\subsection*{Data availability statement}

Not applicable.

\subsection*{Code availability}

Not applicable.

\subsection*{Consent to participate}

Not applicable.

\subsection*{Consent for publication}

Not applicable.

\end{document}